\documentclass{svjour3}
\smartqed  
\usepackage{bm}
\usepackage{amsmath,amssymb,amsbsy,latexsym}
\usepackage{graphicx}
\usepackage{mathptmx}      
\usepackage{epstopdf}
\usepackage{subfigure}
\usepackage{natbib,lineno}
\usepackage{xcolor}
\journalname{BIT}
\date{\today}
\numberwithin{equation}{section}
\numberwithin{figure}{section}
\numberwithin{table}{section}
\newenvironment{rev}{}{}%
\newcommand{\dd}{{\rm d}}

\newcommand\ddh[1]{\tfrac{\dd^{#1}}{\dd h^{#1}}}

\newcommand\nL{\mathcal{L}}
\newcommand\cnL{\check{\mathcal{L}}}
\newcommand\tnL{\tilde{\mathcal{L}}}

\newcommand\nS{\mathcal{S}}
\newcommand\bnS{\bar{\mathcal{S}}}
\newcommand\cnS{\check{\mathcal{S}}}
\newcommand\tnS{\tilde{\mathcal{S}}}

\newcommand\Order{{\mathcal{O}}}

\newcommand\NN{\mathbb N}
\newcommand\method[1]{\mbox{\tt `#1'}}
\newcommand\phh{\phantom{-1}}

\begin{document}

\sloppy
\title{Practical splitting methods for the adaptive integration of nonlinear evolution equations.\\
       Part~I: Construction of optimized schemes and  pairs of schemes}
\titlerunning{Practical splitting methods, Part~I}
\author{Winfried Auzinger \and Harald Hofst{\"a}tter \and David Ketcheson \and Othmar Koch}
\authorrunning{W.\ Auzinger et.\ al.}
\institute{
Winfried~Auzinger \and Harald Hofst{\"a}tter \at
Institut f{\"u}r Analysis und Scientific Computing, Technische Universit{\"a}t Wien \\
Wiedner Hauptstrasse 8--10/E101, 1040 Wien, Austria \\
\email{w.auzinger@tuwien.ac.at, hofi@harald-hofstaetter.at}
\and
David Ketcheson \at
Division of Computer, Electrical and Mathematical Sciences, 4700 King Abdullah University of
Science and Technology (KAUST), Thuwal, 23955-6900, Kingdom of Saudi Arabia \\
\email{David.Ketcheson@kaust.edu.sa}
\and
Othmar~Koch \at
Fakult{\"a}t f{\"u}r Mathematik, Universit{\"a}t Wien \\
Oskar-Morgenstern-Platz 1, 1090 Wien, Austria \\
\email{othmar@othmar-koch.org}
}

\maketitle

\begin{abstract}
We present a number of new contributions to the topic of constructing efficient higher-order
splitting methods for the numerical integration of evolution equations. Particular schemes
are constructed via setup and solution of polynomial systems for the splitting coefficients.
To this end we use and modify a recent approach for generating these systems for
a large class of splittings. In particular, various types of pairs of schemes
intended for use in adaptive integrators are constructed.
\end{abstract}

\section{Introduction} \label{sec:intro}
Operator splitting techniques for the efficient numerical integration of evolution equations
\begin{equation} \label{eq:evolution-equation}
\partial_t u(t) = F(u(t)), \quad t \geq 0, \qquad u(0)~\;\text{given,}
\end{equation}
have become increasingly popular in recent years. Splitting the right-hand side $ F(u) $ into
two or more components in an appropriate way enables efficient and accurate approximations.
In particular, a number of higher-order schemes with real or complex coefficients
have been constructed and analyzed. Relevant contributions to this fields can, e.g., be found in
\cite{blanesetal08,blanesetal13,blanesmoan02,Castella09,chambers03,haireretal02,
mclac95,macqui02,omelyanetal02,yoshida90}.
\begin{rev}Furthermore, application to particular problem classes have been studied in the literature
where the the vector field $ F $ has special properties, such that splitting methods
can be tuned for such cases. In~\cite{blanesetal13a} and~\cite{mclachlan95}, for instance,
perturbations of integrable systems have been considered, say $ F(u) = A(u) + \varepsilon\,B(u) $
where $ \varepsilon $ is s small perturbation parameter. Exploiting this perturbation structure
allows the construction of more efficient (de  facto) higher-order schemes compared to generic ones.
\end{rev}

\paragraph{Overview.}
We present some new contributions to the topic of splitting methods;
\begin{rev}here we will concentrate on the generic case, i.e.,
no special properties of the vector field $ F $ are assumed.\end{rev}
At first we review the approach from~\cite{auzingeretal13c}
for the automatic setup of order conditions
represented by polynomial equations in the coefficients to be determined.
Special cases involving symmetries or composition methods based on lower-order schemes
can be treated as well. Splitting of the right-hand side of~\eqref{eq:evolution-equation}
into two or three components is considered.

The goal is to identify good schemes of a desired order $ p $.
`Good' refers to a compromise between efficiency (minimizing effort) as well as accuracy (minimizing a measure for the
the expected behavior of the local error).
In particular, we focus on the constructions of pairs of schemes
of orders $ (p,p+1) $, where a scheme of order~$ p $ acts as a `worker',
while a related scheme of order $ p+1 $ plays the role of a 'controller' for the
purpose of practical local error estimation.
The idea of using pairs of embedded schemes
(an idea related to Runge-Kutta pairs) is due to~\cite{knth10b}.
Via more flexible embeddings, optimized variants can be constructed.
\begin{rev}Here, `optimization' means searching for schemes where a reasonable measure for
the behavior of the local error becomes minimal among a set of comparable schemes.
It is well-known that this is a very relevant point, because such local error
measures may vary over several orders of magnitude.\end{rev}
We also consider alternative ways of choosing $ (p,p+1) $-pairs,
e.g., adjoint pairs.

\begin{rev}Concerning the search for optimal solutions for a given set of order conditions
(see Section~\ref{sec:OC-solve}), different techniques were applied,
depending on the particular case at hand, including exact, symbolic solution representations
using\footnote{Maple is a product of $ \text{Maplesoft}^{\text{TM}} $.} Maple
(for lower-order schemes), or numerical searches using optimization tools or
straightforward Monte-Carlo techniques.\end{rev}

The ultimate purpose is adaptive integration of evolution equations
based on a reliable local error control.
This topic has been studied in detail, in particular in the context of Schr{\"od}inger equations,
in~\cite{auzingeretal13b,auzingeretal12a,auzingeretal13a,auzingeretal14a}.
In these papers, an alternative method for local error estimation has been constructed
and analyzed. It is based on a computable high order approximation of an integral
representation of the local error in terms of the defect of the numerical solution.
While this approach is rather universal and useful in several cases,
the alternative of using optimized pairs of schemes, if applicable,
will usually be more efficient.

In Part~II of this work we will present a detailed study of adaptive integration,
using both approaches for local error estimation,
for different types of linear and nonlinear evolution equations.

\begin{rev}
\begin{remark}
Recently we became aware of the paper~\cite{blanesetal13a}, where a method of deriving order
conditions has been proposed which is similar to our approach.
Both approaches are based on the notion of a Lyndon basis
(also called Lyndon-Shirshov basis) in a free Lie algebra. In view of the similarities
between our work and~\cite{blanesetal13a}, we stress that we have implemented a fully automatic
computational procedure for deriving order conditions which requires no extra analytical
hand work. This is a versatile implementation, and it can easily be adapted to cover
special cases like palindromic schemes, flexible embeddings, and also splitting into
more than two operators (see Sections~\ref{sec:splitopt} and~\ref{sec:pairs}).

The procedure for setting up higher order conditions
involves the generation of long weighted sums of power products of noncommuting variables
representing the components of the split vector fields.
These sums can easily be distributed in order to obtain a significant speed-up in a parallel environment,
and we have realized such a version.
\end{remark}
\end{rev}

\paragraph{Problem setting and notation.}
For an evolution equation~\eqref{eq:evolution-equation} where the right-hand side
is split into two components,
\begin{equation} \label{eq:AB-problem}
\partial_t u(t) = F(u(t)) = A(u(t)) + B(u(t)), \quad t \geq 0, 
\end{equation}
a single step of a multiplicative splitting scheme, starting from $ u $ and over a step of length $ h $,
is given
by\footnote{$ \phi_F $ denotes the flow associated with the given evolution equation.}
\begin{subequations} \label{eq:AB-scheme}
\begin{equation} \label{eq:AB-scheme-1}
\nS(h,u) = \nS_s(h,\nS_{s-1}(h,\ldots,\nS_1(h,u))) \approx \phi_F(h,u)\,,
\end{equation}
with
\begin{equation} \label{eq:AB-scheme-2}
\nS_j(h,v) = \phi_B(b_j\,h,\phi_A(a_j\,h,v))\,,
\end{equation}
\end{subequations}
with appropriate  coefficients $ a_j,b_j $.
More general schemes based on splitting into three operators are also considered,
see Section~\ref{sec:OC-ABC}, and a special case of additive splitting is also included,
see Section~\ref{sec:OC-special}.

The local error of a splitting step is denoted by
\begin{equation} \label{eq:local_error_notation}
\nS(h,u) - \phi_F(h,u) =: \nL(h,u)\,,
\end{equation}

\paragraph{Contents.}
In Sections~\ref{sec:splitopt} and~\ref{sec:pairs} we describe our approach for
setting up the order conditions for different types of [pairs of] schemes.
Some technical details concerning implementation of this setup procedure
are given in Section~\ref{sec:OC-solve}.
By solving the resulting polynomial systems we have constructed a number of new variants,
and we have compiled a collection of practically relevant
(old and new) schemes and pairs of schemes up to order~$ p=6 $.
This collection can be found at
\begin{flushleft}
\qquad \texttt{http://www.asc.tuwien.ac.at/\~{}winfried/splitting}
\end{flushleft}
and is also expected to be extended in the future,
depending on further investigations on the topic at hand.
We will refer to this webpage throughout as reference~\cite{splithp}
to avoid listing coefficients in the present paper for the sake of brevity.
Some remarks on the schemes collected in~\cite{splithp} are given in Section~\ref{sec:splithp};
\begin{rev}for more detailed information about the properties of the various schemes
we also refer to~~\cite{splithp}.\end{rev}
In Section~\ref{sec:num} we present a numerical example.

\section{Order conditions} \label{sec:splitopt}
Many authors have contributed to the topic of finding good methods.
For an overview on the topic see~\cite{blanesetal08},\,\cite{macqui02}.
Here we do not attempt to describe the relevant approaches and results in detail but
mainly refer to work related to our present activity.
For the relevant mathematical background we refer
to~\cite{blanesetal08,haireretal02,macqui02}.

Among many others, \cite{blanesetal13,blanesmoan02,Castella09},
and~\cite{chambers03} are devoted to the construction of optimal higher-order methods
with real or complex coefficients, either via composition or by solving a set of order conditions
generated in different ways. Order conditions take the form
of a polynomial system in the unknown coefficients
or composition weights $ \omega_\mu $, see Section~\ref{sec:OC-special}.
In the following we recapitulate and illustrate by examples how order conditions
can be set up according to~\cite{auzingeretal13c}; \begin{rev}as mentioned before, this is similar
to one of the approaches taken in~\cite{blanesetal13a}.\end{rev}
Later on we will also present
optimized schemes and pairs of schemes obtained on the basis of this approach,
where `optimized' means that a measure for the local error is chosen as small as possible.

\subsection{Setup of order conditions} \label{sec:setupOC}
There are different ways to generate a polynomial system representing the conditions
on the splitting coefficients for a desired order $ p $.
An essential theoretical basis is the well-known
Baker-Campbell-Hausdorff (BCH) formula, see for example~\cite{haireretal02}.

The approach proposed in~\cite{auzingeretal13c}, which we follow here,
also relies on the BCH formula, but order conditions are set up in a completely automatic way.
Most of the schemes and pairs of schemes specified in~\cite{splithp} have been obtained on the basis of the
algorithm from~\cite{auzingeretal13c}. In the following  we explain
and illustrate this approach by means of examples.
For the purpose of generating order conditions it is sufficient to consider the case of a
linear operator split into two parts $ A $ and $ B $. We denote
\begin{equation*}
A_j = a_j\,A,~~ B_j = b_j\,B, \quad j=1 \ldots s.
\end{equation*}
For the linear case the local error~\eqref{eq:local_error_notation} is of the form
$ \nL(h)\,u $ with a linear operator $ \nL(h) $.

Consider the Taylor expansion of the local
error\footnote{By construction, $ \nL(0)=0 $ for any consistent scheme.}
of a one-step method starting at $ u $,
\begin{equation} \label{eq:lerr-lead}
\nL(h)\,u =
\sum_{q=1}^{p} \tfrac{h^q}{q!}\,\ddh{q}\,\nL(0)\,u
+ \tfrac{h^{p+1}}{(p+1)!}\,\ddh{p+1}\,\nL(0)\,u + \Order(h^{p+2}).
\end{equation}
The method is of order $ p $ iff $ \nL(h) = \Order(h^{p+1}) $; thus the conditions for
order $ p $ are given by
\begin{equation} \label{OC-Tay}
\ddh{}\,\nL(0) = \ldots = \ddh{p}\,\nL(0) = 0.
\end{equation}
For the case of a splitting method we have
(with $ {\bm k}=(k_1,\ldots,k_s) \in \NN_0^s $)
\begin{equation} \label{OC-Tay-1-AB}
\ddh{q}\,\nL(0)
= \sum_{|{\bm k}|=q} \dbinom{q}{{\bm k}}
  \begin{rev}\prod\limits_{j=s \ldots 1}\end{rev}\;
  \sum_{l=0}^{k_j} \dbinom{k_j}{l}\,B_j^{l}\,A_j^{k_j-l}
  \;-\; (A+B)^q,
\end{equation}
If the conditions~\eqref{OC-Tay} are satisfied up to a given order $ p $,
then the leading term of the local error is given by
$ \frac{h^{p+1}}{(p+1)!}\,\ddh{p+1}\,\nL(0) $.
This leading error term is a linear combination of higher-order commutators
of the operators $ A $ and $ B $. As explained in~\cite{auzingeretal13c},
a non-redundant set of order conditions can be built in a recursive way by generating
the symbolic expressions~\eqref{OC-Tay-1-AB} for $ q=1,2,3,\ldots $
in terms of formally linear but non-commuting operators $ A,B $,
and identifying coefficients associated with
power products of $ A $- and $ B $-factors which uniquely identify commutators out of
an appropriate basis of Lie-elements. For this purpose we use the so-called Lyndon basis,
also called Lyndon-Shirshov basis, of the free Lie algebra generated by $ A $ and $ B $.
The elements of this basis are represented by the (associative) Lyndon words
over the alphabet $ \{\mathtt{A,B}\} $, see Table~\ref{TAB:Lyndon-AB}.

\begin{table}[!ht]
\begin{tabular}{|r|r|l|}
\hline
\,$q$\,&$ \ell_q$  &\;Lyndon words over the alphabet $ \{\mathtt{A,B}\} $ $\vphantom{\sum_A^A}$  \\ \hline
  1    & 2      & \;$ {\mathtt{A}},\, {\mathtt{B}} $          $\vphantom{\sum_A^{A^A}}$  \\
  2    & 1      & \;$ {\mathtt{AB}} $              $\vphantom{\sum_A^A}$  \\
  3    & 2      & \;$ {\mathtt{AAB}}, \,{\mathtt{ABB}} $  $\vphantom{\sum_A^A}$  \\
  4    & 3      & \;$ {\mathtt{AAAB}},\, {\mathtt{AABB}},\, {\mathtt{ABBB}} $  $\vphantom{\sum_A^A}$  \\
  5    & 6      & \;$ {\mathtt{AAAAB}},\, {\mathtt{AAABB}},\, {\mathtt{AABAB}},\,
                      {\mathtt{AABBB}},\, {\mathtt{ABABB}},\, {\mathtt{ABBBB}} $  $\vphantom{\sum_A^A}$  \\
  6    & 9      & \;$ {\mathtt{AAAAAB}},\, {\mathtt{AAAABB}},\, {\mathtt{AAABAB}},\,
                      {\mathtt{AAABBB}},\, {\mathtt{AABABB}},\,
                      {\mathtt{AABBAB}},\, {\mathtt{AABBBB}},\, {\mathtt{ABABBB}},\, {\mathtt{ABBBBB}} $ $\vphantom{\sum_A^A}$ \\
  7    & 18     & ~\,\ldots  $\vphantom{\sum_A^A}$  \\
  8    & 30     & ~\,\ldots  $\vphantom{\sum_A^A}$  \\
  9    & 56     & ~\,\ldots  $\vphantom{\sum_A^A}$  \\
 10    & 99     & ~\,\ldots  $\vphantom{\int_{\int_\sum}^A}$  \\
\hline
\end{tabular}
\caption{$ \ell_q $ is the number of words of length~$ q $. \label{TAB:Lyndon-AB}}
\end{table}

Let us first illustrate the procedure by means of a simple example.

\begin{example} \label{example-1}
For $ s=2 $ we have
\begin{subequations}
\begin{align}
\ddh{} \nL(0) &= \underline{(a_1+a_2-1)}\,A + \underline{(b_1+b_2-1)}\,B\,, \label {OC-s=2-1} \\
\ddh{2}\,\nL(0) &= ((a_1+a_2)^2-1)\,A^2  \label {OC-s=2-2} \\
& \quad {} + \underline{(2\,a_2\,b_1-1)}\,AB + (2\,a_1\,b_1 + 2\,a_1\,b_2 + 2\,a_2\,b_2 - 1)\,BA \notag \\
& \quad {} +  ((b_1+b_2)^2-1)\,B^2\,. \notag
\end{align}
\end{subequations}
The basic consistency condition for order $ p=1 $ is $ \ddh{} \nL(0)=0 $ which is
equivalent to $ a_1+a_2=1 $ and $ b_1+b_2=1 $. Assuming these first-order conditions are satisfied,
the second derivative $ \ddh{2}\,\nL(0) $, which now represents the leading error term,
simplifies to the commutator expression
\begin{equation} \label{let-p=2}
\ddh{2}\,\nL(0) = \underline{(2\,a_2\,b_1 - 1)}\,[A,B]\,,
\end{equation}
giving the additional condition $ 2\,a_2\,b_1 = 1$ for order $ p=2 $.
Assuming now that the conditions for $ p=2 $ are satisfied,
the third derivative $ \ddh{3}\,\nL(0) $, which will now represent the leading error term,
is a linear combination of the commutators
$ [A,[A,B]] $ and $ [[A,B],B] $, namely
\begin{equation} \label{let-p=3}
\ddh{3}\,\nL(0) = \underline{(3\,a_2^2\,b_1 - 1)}\,[A,[A,B]]
                   + \underline{(3\,a_2\,b_1^2 - 1)}\,[[A,B],B]\,.
\end{equation}
This computation can be automatized:
\begin{itemize}
\item Generate the representation~\eqref{OC-s=2-1} of $ \ddh{} \nL(0) $
      and extract coefficients of the Lyndon
      words $ \mathtt{A} $ and $ \mathtt{B} $.
      This gives the first-order conditions $ a_1+a_2=1 $ and $ b_1 + b_2= 1 $.
\item Generate the representation~\eqref{OC-s=2-2} of $ \ddh{2} \nL(0) $.
      For a solution of the equations for order~1, the leading local error will have the form
      $ \frac{h^2}{2}\,\ddh{2}\,\nL(0) $ with $ \ddh{2}\,\nL(0) $ from~\eqref{let-p=2}.
      The coefficient of $ [A,B] $ in~\eqref{let-p=2} is determined by extracting
      the coefficient of the Lyndon word $ \mathtt{AB} $ in~\eqref{OC-s=2-2}.
      This gives the equation $ 2\,a_2\,b_1 = 1 $ which,
      together with the first-order conditions, represents a set of conditions for order~$ p=2 $.
\item Generate the representation of $ \ddh{3} \nL(0) $ (we do not display it here).
      For a solution of the equations for order~2, the leading local error will have the form
      $ \frac{h^3}{6}\,\ddh{3}\,\nL(0) $ with $ \ddh{3}\,\nL(0) $ from~\eqref{let-p=3}.
      The coefficients of $ [A,[A,B]] $ and $ [[A,B],B] $ in~\eqref{let-p=3}
      are determined by extracting the coefficients
      of the Lyndon words $ \mathtt{AAB} $ and $ \mathtt{ABB} $
      in the expression for $ \ddh{3} \nL(0) $.
\end{itemize}
In the simple case considered here, there is a one-dimensional manifold of solutions for order~$ p=2 $,
and for each solution $ \{a_1,a_2,b_1,b_2\} $ the size of the coefficients in~\eqref{let-p=3}
is a quality measure.

If a scheme of order~3 is desired, the system of equations is augmented
by the further equations $ 3\,a_2^2\,b_1 = 1 $ and $ 3\,a_2\,b_1^2 = 1 $.
(For the case $ s=2 $ displayed here, the resulting system of equations has no solution;
we need $ s \geq 3 $.)
\end{example}

In general, for arbitrary $ s $ and $ p $, this procedure is continued up to the desired order,
by `implicit recursive elimination' as described in~\cite{auzingeretal13c},
automatically producing a generically non-redundant set of order conditions for a desired order~$ p $.
This process is based on a special bijection between (associative) Lyndon words and bracketed,
non-associative versions of these words which, in our context,
are identified with higher-order commutators representing basis elements for the free Lie algebra
generated by $ A $ and $ B $. The expanded version of such a commutator is a Lie polynomial
in terms of the non-commutative variables $ A $ and $ B $. The essential point is that
{\em its leading monomial, with respect to (alphabetically increasing) lexicographical order, is precisely the
monomial represented by the corresponding Lyndon word;} see~\cite{bokutetal2006}.

In the following, the relation `{\tt <}' refers to lexicographical order of words over the
alphabet $ \{\mathtt{A,B}\} $.


\begin{example} \label{example-2}
Consider a scheme of order $ p=4 $,
i.e., assume that the conditions up to order $ p=4 $ are satisfied.
Then, $ \ddh{5}\,\nL(0) $ is a linear combination of commutators,
or non-associative words, listed below and represented
by the six Lyndon words of length 5 (see Table~\ref{TAB:Lyndon-AB}),
\begin{equation*}
\mathtt{AAAAB < AAABB < AABAB < AABBB < ABABB < ABBBB}\,.
\end{equation*}
The commutators are bracketed, non-associative versions of these
words,\footnote{The bracketing can be computed using the SageMath function
                \texttt{StandardBracketedLyndonWords},
                see \texttt{www.sagemath.org}.}
\begin{align*}
[A,[A,[A,[A,B]]]]
&= \underline{A^4\,B} - 4\,A^3\,B\,A + 6\,A^2\,B\,A^2 - 4\,A\,B\,A^3 + B\,A^4\,, \\
[A,[A,[[A,B],B]]]
&= \underline{A^3\,B^2} - 2\,\underline{\underline{A^2\,B\,A\,B}} + 4\,A\,B\,A\,B\,A - A\,B^2\,A^2 - A^2\,B^2\,A
   - 2\,B\,A\,B\,A^2 + B^2\,A^3\,, \\
[[A,[A,B]],[A,B]]
&= \underline{A^2\,B\,A\,B} - A^2\,B^2\,A -3\,A\,B\,A^2\,B + 4\,A\,B\,A\,B\,A + 2\,B\,A^3\,B - 3\,B\,A^2\,B\,A \\
&  \quad {} - A\,B^2\,A^2 + B\,A\,B\,A^2\,, \\
[A,[[[A,B],B],B]]
&= \underline{A^2\,B^3} - 3\,\underline{\underline{A\,B\,A\,B^2}} + 3\,A\,B^2\,A\,B - 2\,A\,B^3\,A + 3\,B\,A\,B^2\,A
   - 3\,B^2\,A\,B\,A + B^3\,A^2\,, \\
[[A,B],[[A,B],B]]
&= \underline{A\,B\,A\,B^2} - 3\,A\,B^2\,A\,B + 2\,A\,B^3\,A - B\,A^2\,B^2 + 4\,B\,A\,B\,A\,B -3\,B\,A\,B^2\,A \\
&  \quad{} - B^2\,A^2\,B + B^2\,A\,B\,A\,, \\
[[[[A,B],B],B],B]
&= \underline{A\,B^4} - 4\,B\,A\,B^3 + 6\,B^2\,A\,B^2 - 4\,B^3\,A\,B + B^4\,A\,.
\end{align*}
As mentioned above, the \underline{leading (lowest) monomials} in the expanded commutators,
in the sense of lexicographical order,
correspond to the Lyndon words. Note that \underline{\underline{some of these monomials}}
also occur in lower commutators (`lower' again in the sense of lexicographical ordering).
Let us now denote these six commutators by $ K_k,\, k=1 \ldots 6 $. We a priori know that
$ \ddh{5}\,\nL(0) $ is of the form, with $ \ell_5=6 $,
\begin{equation*}
\ddh{5}\,\nL(0) = \sum_{k=1}^{\ell_5}\,\kappa_k\,K_k
\end{equation*}
where the scalars $ \kappa_k $ are multivariate polynomials of degree~5
in the coefficients $ a_j,b_j $ of the underlying scheme of order $ p=4 $.
Therefore the additional conditions for order $ p=5 $ are given by
\begin{subequations}
\label{eq:lamell0}
\begin{equation}
\label{eq:lam0}
\kappa_k = 0, \quad k=1 \ldots \ell_5\,.
\end{equation}
Extracting these coefficients $ \kappa_k $ from the expression~\eqref{OC-Tay-1-AB}
for $ \ddh{5}\,\nL(0) $ is a combinatorial challenge, but we can do better:
We simply extract the coefficients of the Lyndon monomials
\,--\, let us denote them by $ \lambda_k $ \,--\,
which is a standard operation in computer algebra. Now, instead of~\eqref{eq:lam0} we require
\begin{equation}
\label{eq:ell0}
\lambda_k = 0, \quad k=1 \ldots \ell_5\,.
\end{equation}
In our example, for $ {\bm\kappa} = (\kappa_1,\ldots,\kappa_6)^T $
and $ {\bm\lambda} = (\lambda_1,\ldots,\lambda_6)^T $ we have
\begin{equation}
\label{eq:lambdaelleq}
{\bm\lambda} = M\,{\bm\kappa}, \quad \text{with} \quad
M = \left\lgroup\begin{array}{rrrrrr}
      1 & \phh & \phh & \phh & \phh & \phh \\
        &   1  &      &      &      &      \\
        &  -2  &  1   &      &      &      \\
        &      &      &   1  &      &      \\
        &      &      &  -3  &   1  &      \\
        &      &      &      &      &  1
    \end{array}\right\rgroup,
\end{equation}
\end{subequations}
where the lower diagonal entries correspond to the
additional occurrence
of the $ \lambda_k $ in non-leading positions. Therefore the systems~\eqref{eq:lam0}
and~\eqref{eq:ell0} are equivalent.
\end{example}

The situation displayed in this example occurs also in the general case. For any order $ p $,
the vectors $ {\bm\kappa} $ and $ {\bm\lambda} $ consisting of polynomials of degree $ p+1 $
satisfy $ {\bm\lambda} = M\,{\bm\kappa} $
where $ M $ is a lower triangular matrix with unit diagonal. In particular, a Lyndon monomial
$ \lambda_k $ never occurs in an expanded commutator $ K_j $ for $ j>k $ because this would
contradict the leading position~\cite{bokutetal2006} of the Lyndon monomial $ \lambda_j>\lambda_k $ in $ K_j $.


\subsection{Special cases; symmetries} \label{sec:OC-special}
In the sequel,
\begin{equation*}
\nS^\ast(h,u) = \nS^{-1}(-h,u)
\end{equation*}
denotes the adjoint scheme associated with $ \nS $.

The order conditions generated by the algorithm indicated in Section~\ref{sec:setupOC}
are generically non-redundant. However, there exist special cases:
\begin{itemize}
\item {\em Symmetric}\, (or: `time-symmetric') one-step schemes are characterized by the property
      \begin{equation} \label{eq:symmetric-relation}
      \nS(-h,\nS(h,u)) = u, \quad \text{i.e.,} \quad \nS(h,u) = \nS^\ast(h,u).
      \end{equation}
      For symmetric splitting schemes we have
      either $ a_1=0 $ or $ b_s=0 $, and the remaining coefficient tupels
      $ (a_j) $ and $ (b_j) $ are both palindromic.
      Since symmetric schemes have an even order $ p $
      (cf.~\cite[Chapter~3]{haireretal02}),
      only odd-order conditions for an appropriately reduced number of
      free coefficients need to be imposed. The general algorithm described in
      Section~\ref{sec:setupOC} can easily be adapted to this case.
\item The following type of schemes seems not to have been considered earlier in the
      literature:\footnote{The Lie-Trotter scheme, with $ s=p=1,\; a_1=b_1=1 $,
                           is a trivial special case.}

      {\em Palindromic}\, schemes, or `reflected schemes' in the terminology of~\cite{auzingeretal13c},
      are characterized by $ b_j = a_{s+1-j},\, j=1 \ldots s $, i.e.,
      \begin{equation} \label{eq:palicoe}
      \begin{aligned}
        & (a_1,b_1,a_2,b_2,\ldots,a_{s-1},b_{s-1},a_s,b_s) \\
      =~& (a_1,b_1,a_2,b_2,\ldots,b_2,~~~\,a_2,~~~\,b_1,a_1).
      \end{aligned}
      \end{equation}
      Assume a scheme of order $ p $ is given, and consider a splitting
      step of the form~\eqref{eq:AB-scheme}. Interchanging the roles of $ A $ and $ B $,
      i.e., replacing \eqref{eq:AB-scheme} by
      \begin{subequations} \label{eq:BA-scheme}
      \begin{equation} \label{eq:BA-scheme-1}
      \cnS(h,u) = \cnS_s(h,\cnS_{s-1}(h,\ldots,\cnS_1(h,u))),
      \end{equation}
      with
      \begin{equation} \label{eq:BA-scheme-2}
      \cnS_j(h,v) = \phi_A(b_j\,h,\phi_B(a_j\,h,v)),
      \end{equation}
      \end{subequations}
      also results in a scheme of order $ p $.
      If $ \nS $ is palindromic then
      \begin{equation} \label{eq:palindromic-relation}
      \nS(-h,\cnS(h,u)) = u, \quad \text{i.e.,} \quad \cnS(h,u) = \nS^\ast(h,u).
      \end{equation}
      Thus we infer from~\cite[Theorem~II.3.2]{haireretal02}
      that in the palindromic case the local errors
      $ \nL(h,u) = \nS(h,u) - \phi_F(h,u) $ and $ \cnL(h,u) = \cnS(h,u) - \phi_F(h,u) $
      are related via
      \begin{subequations} \label{eq:le-pali}
      \begin{align}
      \nL(h,u)  &= ~~~~~~~~~~\,C(u)\,h^{p+1} + \Order(h^{p+2}), \label{eq:le-pali1} \\
      \cnL(h,u) &= (-1)^p\,C(u)\,h^{p+1} + \Order(h^{p+2}), \label{eq:le-pali2}
      \end{align}
      \end{subequations}
      with $ C(u) = \tfrac{1}{(p+1)!}\,\ddh{p+1}\,\nL(0,u) $.
      For an ansatz with palindromic coefficients,
      exchanging the roles of $ A $ and $ B $ in the algorithm from
      Section~\ref{sec:setupOC} will lead to the identical set of order conditions.
      Therefore the order conditions associated with `Lyndon twins'
      are pairwise identical. Here, we call a pair of Lyndon words a twin
      if one of them is obtained by exchanging the role of $ A $ and $ B $
      and reading it from right to
      See Table~\ref{TAB:Lyndon-AB};
      for instance, the 6~words of odd length~5 consist of three twins;
      the 9~words of even length~6 consist of three twins, the selfie
      $ \mathtt{AAABBB} $, and two solitary words.

      Due to this redundancy the number of order conditions is appropriately reduced.
\item Higher order one-step schemes
      can be generated by $ m $-fold {\em composition}\, of lower-order schemes with
      appropriately chosen sub-steps $ h_\mu = \omega_\mu h $ satisfying
      $ \omega_1 + \ldots + \omega_m = 1 $
      plus additional conditions guaranteeing that a certain order is
      obtained.\footnote{
        We note that the idea of composition is of a general nature
        and not restricted to the class of splitting methods.}

      A popular class of composition methods
      are symmetric Strang compositions. Schemes of this type of orders $ 4,6 $ and higher
      were first devised in~\cite{yoshida90}.
      Some of the composition  coefficients have to be chosen negative, and the local error
      measures of these composition schemes are rather large. On the other hand,
      for higher orders, composition beats the generic lower limits on the number $ s $
      of stages such that a given order $ p $ can be expected. For instance,
      the 7-fold 6-th order symmetric Strang composition \cite[`Y 8-6']{splithp} recombines into
      an 8-stage scheme, whereas the generic number of order conditions for a symmetric scheme
      of order~$ p=6 $ is~10, which would require $ s=10 $ stages involving 11~free coefficients.

      \begin{rev}
      Evidently, (symmetric) compositions are an attractive option for
      constructing higher-order schemes. Therefore we have included this class into
      our considerations concerning the search for optimal variants (see Section~\ref{sec:OC-solve}).
      \end{rev}
\end{itemize}

\subsection{Complex coefficients} \label{sec:OC-complex}
Our considerations are not restricted to schemes with real coefficients $ a_j,b_j $.
Complex schemes, with coefficients having positive real parts, are appropriate for the application
of splitting methods to parabolic problems, since real schemes with positive coefficients
do not exist for order $ p \geq 3 $, see~\cite{blanesetal08}.
For this class of methods, in particular based on complex compositions,
we refer to~\cite{chambers03} and~\cite{blanesetal13}.

\subsection{Splitting into more than two operators} \label{sec:OC-ABC}
We also consider evolution equations where the right-hand side splits into three parts,
\begin{equation} \label{eq:ABC-problem}
\partial_t u(t) = F(u(t)) = A(u(t)) + B(u(t)) + C(u(t)), \quad t \geq 0, \qquad u(0)~\text{given,}
\end{equation}
and according multiplicative splitting schemes,
\begin{subequations} \label{eq:ABC-scheme}
\begin{equation} \label{eq:ABC-scheme-1}
\nS(h,u) = \nS_s(h,\nS_{s-1}(h,\ldots,\nS_1(h,u))) \approx \phi_F(h,u),
\end{equation}
with
\begin{equation} \label{eq:ABC-scheme-2}
\nS_j(h,v) = \phi_C(c_j\,h,\phi_B(b_j\,h,\phi_A(a_j\,h,v))).
\end{equation}
\end{subequations}
The methodology from~\cite{auzingeretal13c} can be directly generalized to the case of
splitting into more than two operators. For the practically relevant case of splitting into
three operators $ A,B,C $, as in~\eqref{eq:ABC-scheme},
the representation~\eqref{OC-Tay-1-AB} generalizes as follows,
with $A_j = a_j\,A,\, B_j=b_j\,B,\, C_j=c_j\,C $, and
$ {\bm k}=(k_1,\ldots,k_s) \in \NN_0^s $,\, $ {\bm l} = (l_A,l_B,l_C) \in \NN_0^3 $:
\begin{equation} \label{OC-Tay-1-ABC}
\ddh{q}\,\nL(0)
= \sum_{|{\bm k}|=q} \dbinom{q}{{\bm k}}
  \begin{rev}\prod\limits_{j=s \ldots 1}\end{rev}\;
  \sum_{|{\bm l}|=k_j} \dbinom{k_j}{{\bm l}}\,C_j^{l_C}\,B_j^{l_B}\,A_j^{l_A}
  \;-\; (A+B+C)^q.
\end{equation}
On the basis of these identities, the algorithm from Section~\ref{sec:setupOC} generalizes in a straightforward way.
The Lyndon basis representing independent commutators now corresponds to Lyndon words
over the alphabet $ \{\mathtt{A,B,C}\} $, see Table~\ref{TAB:Lyndon-ABC}.

Concerning symmetries, similar considerations as in Section~\ref{sec:OC-special} apply.

\begin{table}[!ht]
\begin{tabular}{|r|r|l|}
\hline
\,$q$\,&$ \ell_q$  &\;Lyndon words over the alphabet $ \{\mathtt{A,B,C}\} $ $\vphantom{\sum_A^A}$  \\ \hline
  1    & 3      & \;$ {\mathtt{A}},\, {\mathtt{B}},\, {\mathtt{C}} $  $\vphantom{\sum_A^{A^A}}$  \\
  2    & 3      & \;$ {\mathtt{AB}},\, {\mathtt{AC}},\, {\mathtt{BC}} $              $\vphantom{\sum_A^A}$  \\
  3    & 8      & \;$ {\mathtt{AAB}},\, {\mathtt{AAC}},\,
                      {\mathtt{ABB}},\, {\mathtt{ABC}},\,
                      {\mathtt{ACB}}, \,{\mathtt{ACC}},\,
                      {\mathtt{BBC}}, \,{\mathtt{BCC}} $ $\vphantom{\sum_A^A}$  \\
  4    & 18     & ~\,\ldots  $\vphantom{\sum_A^A}$  \\
  5    & 48     & ~\,\ldots  $\vphantom{\sum_A^A}$  \\
  6    & 115    & ~\,\ldots  $\vphantom{\sum_A^A}$  \\
  7    & 312    & ~\,\ldots  $\vphantom{\sum_A^A}$  \\
  8    & 810    & ~\,\ldots  $\vphantom{\sum_A^A}$  $\vphantom{\int_{\int_\sum}^A}$ \\
\hline
\end{tabular}
\caption{$ \ell_q $ is the number of words of length~$ q $. \label{TAB:Lyndon-ABC}}
\end{table}

For a general convergence theory of ABC-splitting
for the linear case and some applications we refer to~\cite{auzingeretal14a}.
For example, splitting into three operators can be used to handle evolution equations
where the right-hand side splits up into two non-autonomous parts. Introducing the independent
variable~$ t $ as an unknown variable satisfying $ t'=1 $, such a problem can be
formally considered as an autonomous system split into three parts.
In this case, splitting means that the variable $ t $ is frozen over several subintervals
comprising an integration step.
Since the ODE $ t'=1 $ is trivial, a large number of higher-order commutators vanishes
in this case, and therefore the number of necessary order conditions is significantly reduced.
This special situation will be considered in detail later on.

\section{Pairs of splitting schemes} \label{sec:pairs}
For the purpose of efficient local error estimation as a basis for adaptive stepsize selection,
using pairs of related schemes is a well-established idea. One of the schemes, of order $ p $,
acts as the worker, and the other, of order $ p+1 $,
is the controller responsible for local error
estimation.\footnote{Of course, a scheme acting as a controller can also be used as an integrator in a normal way.}
Criteria for the selection of pairs of schemes are accuracy and computational efficiency.

Order conditions for pairs of schemes of the types listed below
can be generated with minor modifications of the approach described in Section~\ref{sec:splitopt}.
\begin{itemize}
\item {\em Embedded pairs.}
In~\cite{knth10b}, pairs of splitting schemes of orders $ p $ and $ p+1 $ are specified.
The idea is to select a controller $ \bnS $ of order $ p+1 $ and to construct a worker $ \nS $
of order~$ p $ for which a maximal number of stages $ \nS_j $ coincides with those of the controller.
Let $ a_j,b_j $ and $ \bar{a}_j,\bar{b}_j $ denote the coefficients of the worker and controller,
respectively. The approach adopted in~\cite{knth10b} may be called static, finding $ \nS $ and $ \bnS $
such that $ a_j={\bar a}_j $ and $ b_j={\bar b}_j $ for as many $ j=1,2,\ldots $ as possible.
In this sense the schemes are related to each other but, in general,
the total number of order conditions, and thus the total number of necessary evaluations,
is the same as for an arbitrary unrelated $ (p,p+1) $ pair.

Here we develop the idea of embedding further: Again we fix a `good' controller of order $ p+1 $
and wish to adjoin to it a `good' worker of order $ p $.
Since the number of stages $ \bar{s} $ of $ \bnS $
will be higher than the number of stages $ s $ of $ \nS $,
we can select an optimal embedded worker
$ \nS $ from a set of candidates obtained by flexible embedding,
where the number of coinciding coefficients is not a priori fixed.
\begin{example}
In~\cite{knth10b}, an embedded $ (3,4) $-pair was constructed, where the controller is an optimized
symmetric scheme of order $ p=4 $ with $ s=7 $ stages due to~\cite{blanesmoan02}, with
local error measure LEM=0.01 (`LEM' in the sense of~\eqref{eq:LEM-3} below).
The worker specified in~\cite{knth10b} is a scheme of order $ p=3 $ with $ s=6 $
stages, where the coefficients $ a_1,a_2,a_3,a_4 $ and $ b_1,b_2,b_3 $ coincide with those
of the controller. This amounts to 7 additional evaluations for the worker, and its local error
measure is LEM=0.2.

For flexible embedding, in contrast, we consider all possible embedded workers, and we find that
a scheme of order $ p=3 $ with $ s=4 $ stages is to be preferred,
see~\cite[Emb 4/3 BM PRK/A]{splithp},
where $ a_1,a_2 $ and $ b_1 $ coincide with those of the controller.
This amounts to 5 additional evaluations for the worker, and the controller has LEM=0.1.
\end{example}
\item {\em Milne pairs.} In the context of multistep methods for ODEs,
the so-called Milne device is a well-established technique for constructing
pairs of schemes. In our context, one may aim for finding a pair
$ (\nS,\tnS) $ of schemes of the same type, with equal $ s $ and $ p $, such that their local errors
$ \nL,\tnL $ are related according to
\begin{subequations} \label{eq:le-milne}
\begin{align}
\nL(h,u)  &= ~~C(u)\,h^{p+1} + \Order(h^{p+2}), \label{eq:le-milne1} \\
\tnL(h,u) &= \gamma\,C(u)\,h^{p+1} + \Order(h^{p+2}), \label{eq:le-milne2}
\end{align}
\end{subequations}
with $ \gamma \not= 1 $. Then, the additive scheme
\begin{equation*}
{\bar\nS}(h,u) = -\tfrac{\gamma}{1-\gamma}\,\nS(h,u) + \tfrac{1}{1-\gamma}\,\tnS(h,u)
\end{equation*}
is a method of order $ p+1 $, and
\begin{equation*}
\nS(h,u) - {\bar\nS}(h,u) = \tfrac{1}{1-\gamma}\,\big( \nS(h,u) - \tnS(h,u) \big)
\end{equation*}
provides an asymptotically correct local error estimate for $ \nS(h,u) $.

\smallskip
\item \begin{rev}{\em Adjoint pairs.}\end{rev}
Let $ \nS $ be a scheme of odd order $ p $ and and $ \nS^\ast $ its adjoint,
see Section~\ref{sec:OC-special}. Due to~\cite[Theorem~II.3.2]{haireretal02}
the leading error terms of $ \nS $ and its adjoint $ \nS^\ast  $
are identical up to the factor $ -1 $. Therefore, the averaged additive scheme
\begin{equation}
\label{eq:pali-avaraged}
{\bar\nS}(h,u) = \tfrac{1}{2}\,\big( \nS(h,u)+\nS^\ast(h,u) \big)
\end{equation}
is a method of order\footnote{For the simplest case of the Lie-Trotter scheme was
   already observed in~\cite{strang68}.} $ p+1 $, and
\begin{equation*}
\nS(h,u) - {\bar\nS}(h,u) = \tfrac{1}{2}\,\big( \nS(h,u) - \nS^\ast(h,u) \big)
\end{equation*}
provides an asymptotically correct local error estimate for $ \nS(h,u) $.
In this case the additional effort for computing the local error estimate is
identical with the effort for the worker $ \nS $ but not higher as is the case
for embedded pairs.
\begin{rev}An example are palindromic pairs, where $ \nS $ is palindromic (of odd order $ p $),
such that $ \nS^\ast = \cnS $, see Section~\ref{sec:OC-special}.
\end{rev}
\end{itemize}
For detailed comments on a number of new pairs listed in~\cite{splithp}, see Section~\ref{sec:splithp}.

\section{Implementation aspects: constructing schemes and minimizing local error terms} \label{sec:OC-solve}
Our approach for setting up order conditions described in Section~\ref{sec:setupOC} has
been implemented in Maple~18. We use the {\tt Physics} package for the manipulation of noncommuting symbols,
and tables of Lyndon words generated using an algorithm devised in~\cite{duval88}.
Since the number of terms in~\eqref{OC-Tay-1-AB} resp.~\eqref{OC-Tay-1-ABC} rapidly increases with $ q $
we have implemented a parallel version relying on Maple's {\tt Grid} package.
In particular, the job of generating all the terms
in the long sums~\eqref{OC-Tay-1-AB} and~\eqref{OC-Tay-1-ABC}
can be (equi-)distributed over several parallel threads.

The resulting set of order conditions is a multivariate polynomial system which, for higher orders,
requires numerical solution techniques.
Once a scheme of order $ p $ has been found, its leading local error term is of the form
(see Section~\ref{sec:splitopt})
\begin{equation} \label{locerr-lem}
\tfrac{h^{p+1}}{(p+1)!}\,\ddh{p+1} \nL(0) = \sum_{k=1}^{\ell_{p+1}} \kappa_{p+1,k}\,K_{p+1,k},
\end{equation}
with $ \ell_{p+1} $ commutators $ K_{p+1,k} $ associated with Lyndon words of length $ p+1 $.
To compare schemes of equal order $ p $ one may consider
\begin{subequations} \label{eq:LEM}
\begin{equation} \label{eq:LEM-2}
\Big(\,\sum_{k=1}^{\ell_{p+1}} |\kappa_{p+1,k}|^2 \Big)^{1/2}
\end{equation}
as a reasonable measure for the accuracy of a scheme. However, we use the quantity
\begin{equation} \label{eq:LEM-3}
\text{LEM} := \Big(\,\sum_{k=1}^{\ell_{p+1}} |\lambda_{p+1,k}|^2 \Big)^{1/2}
\end{equation}
\end{subequations}
instead. Using~\eqref{eq:LEM-3} has the advantage that
the coefficients $ \lambda_k = \lambda_{p+1,k} $ are exactly those
which are generated in the course of the setup
of the conditions for order $ p+1 $, see Section~\ref{sec:setupOC},
while the coefficients from~\eqref{eq:LEM-2} are more difficult to compute
(cf.\ the discussion in Section~\ref{sec:setupOC}).
Since different particular solutions to the order conditions typically result in leading local error terms
varying over several orders of magnitudes, we consider~\eqref{eq:LEM-3} equally reasonable as~\eqref{eq:LEM-2}.

For finding and evaluating solutions and pairs of solutions we follow two different strategies.
\begin{itemize}
\item
For the case where the number of equations equals the number of free coefficients
we expect a set of isolated solutions. In this case we use the {\tt fsolve} function
in Maple combined with a Monte-Carlo strategy for generating different initial intervals.
Higher precision is used to generate solutions with double precision accuracy.
For each detected solution the LEM~\eqref{eq:LEM-3} is computed.
\item
Especially for the case where the number of equations is smaller than the number of free coefficients,
the problem is to be considered as a constrained minimization problem:
Minimize the LEM representing the objective function,
with the order conditions imposed as nonlinear equality constraints.
To this end we employ state-of-the-art techniques
which have also been applied for the construction of special classes of Runge-Kutta methods,
see for instance~\cite{ketchetal13}.
In particular we have used the MATLAB\footnote{MATLAB is a trademark of The\,Math\,Works, Inc.}
optimizer {\tt fmincon}.
Again a large number of initial guesses are generated randomly,
since this optimization problem is nonconvex in general.
The results cannot be guaranteed globally optimal,
but results from an exhaustive search usually suggest that this is indeed the case.

A post-processing, i.e., refining the solutions to full double precision,
is again performed in Maple using higher precision {\tt sfloat} arithmetic.
\end{itemize}
We have also re-checked a number of known methods, refined their coefficients to
full double precision, and computed their LEMs.

\section{Schemes from the collection~\cite{splithp}} \label{sec:splithp}
This collection is not intended to be exhaustive. It includes some known and quite a number of
new schemes, in particular pairs of schemes, up to order~$ p=6 $, with their essential properties.
Some methods are included mainly for the sake of completeness or their historical
significance.

In the following we comment on some of these methods;
for complete information, consult~\cite{splithp}. `Best' or `optimal'
means that it has minimal LEM~\eqref{eq:LEM-3} among a certain class of methods with comparable
effort for a given order $ p $.
In some simple cases such optimality properties can be established theoretically;
for higher orders we have resorted to more or less exhaustive numerical search.

Methods whose label contains the letter {\tt `A'} are new, or taken again into
consideration in the context of constructing pairs,
or their LEM has been computed for the first
time.\footnote{Of course, `new' may not be considered as a rigorous statement in each case since
               the literature on the subject is rather large by now.}
The list also includes some pairs of embedded schemes (\method{Emb ...}),
pairs of Milne type (\method{Milne ...}), and
palindromic pairs (\method{PP ...}), see Section~\ref{sec:pairs}.

\begin{rev}
More detailed information about all these methods can be found on the webpage~\cite{splithp}.
\end{rev}

\begin{rev}
\begin{remark}
In several cases we observed that palindromic schemes tend to have minimal LEMs among a
set of comparable schemes, for instance the third-order scheme in the pair
\method{PP 3/4 A}. This is the reason why we have included some adjoint pairs of
(optimized) palindromic type of orders $ (p,p+1) $ (with $ p $ odd) in our collection.
\end{remark}
\end{rev}

\subsection{Splitting into two operators (`AB schemes')}
\paragraph{Real coefficients.}
\begin{itemize}
\item The best schemes up to order $ p=5 $ we have found are palindromic:
      \begin{itemize}
      \item \method{best 2-stage 2nd order} ($ s=p=2 $).
      \item \method{Emb 3/2 AKS} (palindromic controller with $ s=p=3 $).
      \item \method{Emb 4/3 AKS p} (palindromic controller with $ s=5 $, $ p=4 $).

            In particular, this scheme has essentially the same LEM as the fourth order scheme
            from~\cite{blanesmoan02} which has been used in~\cite{knth10b},
            but it has only~5 stages instead of~7.
      \item \method{Emb 5/4 A} (palindromic controller with $ s=8 $, $ p=5 $),
            see also \method{PP 5/6 A}.
      \end{itemize}
\item \method{Emb 5/4 AK (ii)} is an optimized embedded pair. The controller is a new
      scheme with $ s=7 $, $ p=5 $, and the worker of order $ p=4 $ is chosen out of several
      dozens of candidates of order $ 4 $ which share the same computational effort
      but have LEMs varying over several orders of magnitudes.
\item Palindromic pairs: \method{PP 3/4 A}, \method{PP 5/6 A}.
\end{itemize}

\paragraph{Complex coefficients (with positive real parts).}
\begin{itemize}
\item Since for order $ p=3 $ we need 5~conditions, the question is whether there
      exists a third-order scheme with $ s=3 $ and 5~evaluations. It turns out that
      the only scheme of this type, \method{A 3-3 c}, has complex coefficients.
\item \method{A 4-4 c} ($s=4 $, $ p=4 $)
      is the best complex symmetric Strang composition method of order~$ 4 $;
      see also~\cite{Castella09} and~\cite{chambers03}.
\item \method{Emb 3/2 A c} and \method{Emb 4/3 A c} are embedded pairs with
      palindromic controller and optimized worker.
      We note that the controller in \method{Emb 4/3 A c} ($ s=5 $, $p=4$)
      has a significantly smaller LEM than \method{A 4-4 c} (factor $ \approx 20 $).
\item \method{C 8-6 c} ($ s=8 $, $ p=6 $)
      is the best symmetric complex Strang composition method of order~$ 6 $;
      see also~\cite{Castella09} and~\cite{chambers03}.
\item Palindromic pairs: \method{PP 3/4 A c}, \method{PP 5/6 A c}.
\end{itemize}

\subsection{Splitting into three operators (`ABC schemes')}
Due to the rapidly increasing number of generic order conditions,
finding general higher order schemes would be a very challenging task \begin{rev}for this case.\end{rev}
For $ p=6 $, for instance, the generic number of order conditions is 196 for the general case
and 59 for the symmetric case.
\begin{rev}For $ p=6 $ we therefore only consider real or
complex Strang compositions which are easier to construct and lead to more compact schemes.\end{rev}
Generating the expression for the leading error term $ \ddh{7} \nL(0) $ for the purpose of
computing the LEM for $ p=6 $, involving 312~coefficients
(see Table~\ref{TAB:Lyndon-ABC}), is computationally expensive,
\begin{rev}but it can be done at reasonable effort,
for the purpose of computing the LEM of a given composition and comparing different variants.\end{rev}

\paragraph{Real coefficients.}
\begin{itemize}
\item \method{AK 5-2} ($ s=5 $, $ p=2 $, 9 evaluations) appears to be a possible
      rival of the Strang scheme ($ s=3 $, $ p=2 $, 5 evaluations),
      with a LEM which is smaller by a factor $ \approx 7 $.
\item \method{PP 3/4 A 3} is a palindromic pair based on the best palindromic
      scheme found for $ s=6 $, $ p=3 $.
\item \method{Y 7-4} ($ s=7, p=4 $, 13 evaluations)
      is the best symmetric Strang composition of order~$ p=4 $.
      It is the analog of the AB composition \method{Y 4-4}, with the same composition weights.
\item \method{AK 11-4} ($ s=11 $, $ p=4 $, 21 evaluations) has been found on the basis of 11 conditions
      for a symmetric ABC scheme of order~4. Its LEM is smaller by a factor $ \approx 13 $
      compared to \method{Y 7-4}.
\item \method{AY 15-6} ($ s=15, p=6 $)
      is the best symmetric Strang composition of order~$ p=6 $.
      It is the analog of the AB composition \method{Y 8-6}, with the same composition weights.
\end{itemize}

\paragraph{Complex coefficients (with positive real parts).}
\begin{itemize}
\item \method{AK 7-4 c} ($ s=7, p=4 $)
      is the best symmetric Strang composition of order~$ p=4 $.
      It is the analog of the AB composition \method{A 4-4-c}, with the same composition weights.
\item \method{AK 15-6 c} ($ s=15, p=6 $)
      is the best symmetric Strang composition of order~$ p=6 $.
      It is the analog of the AB composition \method{C 8-6-c}, with the same composition weights.
\end{itemize}

\section{Numerical example} \label{sec:num}
For a numerical illustration, in particular concerning the expected performance of palindromic schemes,
we consider an example of a system of coupled nonlinear evolution equations of Schr{\"o}dinger
type (see~\cite{wadatietal92}),
\begin{equation} \label{eq:schroesys}
\begin{aligned}
\mathrm{i}\,\Big(\frac{\partial \psi_1}{\partial t}
                 + \delta \frac{\partial \psi_1}{\partial x}\Big)
+ \frac{1}{2}\,\frac{\partial^2 \psi_1}{\partial x^2}
+ \big( |\psi_1|^2 + e\,|\psi_2|^2 \big) \psi_1 &= 0, \\[\jot]
\mathrm{i}\,\Big(\frac{\partial \psi_2}{\partial t}
                 - \delta \frac{\partial \psi_2}{\partial x}\Big)
+ \frac{1}{2}\,\frac{\partial^2 \psi_2}{\partial x^2}
+ \big( e\,|\psi_1|^2 + |\psi_2|^2 \big) \psi_2 &= 0,
\end{aligned}
\end{equation}
with initial condition chosen such that the exact solution is a pair of solitons,
\begin{equation*}
\begin{aligned}
\psi_1(x,t) &=
\frac{\sqrt{2\,\beta}}{1+e}\,\text{sech}\big( \sqrt{2\,\beta}\,(x-v\,t) \big)\,
\mathrm{e}^{\mathrm{i}\left((v-\delta)x + (\beta-(v^2-\delta^2)/2)\,t\right)}, \\[\jot]
\psi_2(x,t) &=
\frac{\sqrt{2\,\beta}}{1+e}\,\text{sech}\big( \sqrt{2\,\beta}\,(x-v\,t) \big)\,
\mathrm{e}^{\mathrm{i}\left((v+\delta)x + (\beta-(v^2-\delta^2)/2)\,t\right)},
\end{aligned}
\end{equation*}
which is exponentially decreasing with $ |x| $. We start at $ t=0 $,
the parameters are chosen as $ \delta=0.5 $, $ \beta=1.0 $, $ v=1.1 $, and $ e=0.8 $.

We impose periodic boundary conditions on the interval $ x_{min},x_{max}=[-50,70] $
using an equidistant grid of size 2\,048.
For splitting we choose the time step $ h $ and separately integrate
\begin{itemize}
\item the kinetic part (`A') involving the derivatives
      w.r.t.\ $ x $, using a Fourier spectral discretization,
\item and the nonlinear `ODE part' (`B'), which can be exactly propagated:
      At each grid point~$ x $, the respective solution
      $ (\psi_{1,B},\psi_{2,B}) = (\psi_{1,B}(x,t),\psi_{2,B}(x,t)) $ of the ODE system
      \begin{align*}
      \mathrm{i}\,\frac{\mathrm{d} \psi_{1,B}}{\mathrm{d}\,t}
      + \big( |\psi_{1,B}|^2 + e\,|\psi_{2,B}|^2 \big) \psi_{1,B} &= 0, \\
      \mathrm{i}\,\frac{\mathrm{d} \psi_{2,B}}{\mathrm{d}\,t}
      + \big( e\,|\psi_{1,B}|^2 + |\psi_{2,B}|^2 \big) \psi_{2,B} &= 0,
      \end{align*}
      starting at $ t_0 $ is given by
      \begin{align*}
      \psi_{1,B}(x,t) &=
      \mathrm{e}^{\,\mathrm{i}\,(t-t_0)\,\left(|\psi_{1,B}(x,t_0)|^2 + e\,|\psi_{2,B}(x,t_0)|^2\right)} \psi_{1,B}(x,t_0), \\
      \psi_{2,B}(x,t) &=
      \mathrm{e}^{\,\mathrm{i}\,(t-t_0)\,\left(e\,|\psi_{1,B}(x,t_0)|^2 + |\psi_{2,B}(x,t_0)|^2\right)} \psi_{2,B}(x,t_0).
      \end{align*}
\end{itemize}
All computations were performed in standard double precision arithmetic.
In Tables~\ref{TAB:schroesys-1} and~\ref{TAB:schroesys-2}, `err' refers to a canonically scaled
discrete $ L_2 $\,-\,norm, and `ord' refers to the order observed.


\begin{table}[!ht]
\begin{tabular}{|l||l|l||l|l|}
\hline
 & \multicolumn{2}{|l||}{scheme (i)}
 & \multicolumn{2}{|l|}{scheme ((i)+(ii))/2}
$\vphantom{\int_A^A}$ \\ \hline
$ h $ & $\text{err}_{local}$ & $\text{ord}_{local}$
      & $\text{err}_{local}$ & $\text{ord}_{local}$
$\vphantom{\int_A^A}$  \\ \hline
0.100\,E$+$00 & 0.524\,E$-$03 &      & 0.120\,E$-$03 & $\vphantom{\int^A}$  \\
0.500\,E$-$01 & 0.374\,E$-$04 & 3.74 & 0.467\,E$-$05 & 4.69 \\
0.250\,E$-$01 & 0.246\,E$-$05 & 3.93 & 0.150\,E$-$06 & 4.96 \\
0.125\,E$-$01 & 0.156\,E$-$06 & 3.98 & 0.468\,E$-$08 & 5.01 \\
0.625\,E$-$02 & 0.982\,E$-$08 & 3.99 & 0.146\,E$-$09 & 5.00 \\
0.313\,E$-$02 & 0.614\,E$-$09 & 4.00 & 0.455\,E$-$11 & 5.00 \\
0.156\,E$-$02 & 0.384\,E$-$10 & 4.00 & 0.142\,E$-$12 & 5.00 \\
0.781\,E$-$03 & 0.240\,E$-$11 & 4.00 & 0.456\,E$-$14 & 4.96 \\
\hline
\end{tabular}
\qquad
\begin{tabular}{|l||l|l|}
\hline
\multicolumn{2}{|l|}{scheme (i)}
$\vphantom{\int_A^A}$ \\ \hline
$\text{err}_{global}$ & $\text{ord}_{global}\!\!\!\!\!$
$\vphantom{\int_A^A}$  \\ \hline
0.165\,E$-$02 & $\vphantom{\int^A}$  \\
0.106\,E$-$03 & 3.96 \\
0.912\,E$-$05 & 3.54 \\
0.100\,E$-$05 & 3.18 \\
0.123\,E$-$06 & 3.03 \\
0.154\,E$-$07 & 2.99 \\
0.194\,E$-$08 & 2.99 \\
0.244\,E$-$09 & 2.99 \\
\hline
\end{tabular}
\caption{Error tables for the palindromic pair \method{PP 3/4 A}
         applied to problem~\eqref{eq:schroesys}. \newline
         {\em Left:}\, Local error (first step) for scheme~(i) starting with `A'
         of order~$ 3 $, and for the averaged scheme (see~\eqref{eq:pali-avaraged})
         of order~$ 4 $. \newline
         {\em Right:}\, Global error for scheme~(i) at $ t_{end} = 5.0 $.
         \label{TAB:schroesys-1}}
\end{table}

\begin{table}[!ht]
\begin{tabular}{|l||l|l||l|l|}
\hline
 & \multicolumn{2}{|l||}{scheme (i)}
 & \multicolumn{2}{|l|}{scheme ((i)+(ii))/2}
$\vphantom{\int_A^A}$ \\ \hline
$ h $ & $\text{err}_{local}$ & $\text{ord}_{local}$
      & $\text{err}_{local}$ & $\text{ord}_{local}$
$\vphantom{\int_A^A}$  \\ \hline
0.100\,E$+$00 & 0.322\,E$-$04 &      & 0.318\,E$-$04 & $\vphantom{\int^A}$  \\
0.500\,E$-$01 & 0.590\,E$-$06 & 5.77 & 0.578\,E$-$06 & 5.78 \\
0.250\,E$-$01 & 0.723\,E$-$08 & 6.35 & 0.625\,E$-$08 & 6.53 \\
0.125\,E$-$01 & 0.903\,E$-$10 & 6.32 & 0.534\,E$-$10 & 6.87 \\
0.625\,E$-$02 & 0.129\,E$-$11 & 6.13 & 0.427\,E$-$12 & 6.97 \\
\hline
\end{tabular}
\qquad
\begin{tabular}{|l||l|l|}
\hline
\multicolumn{2}{|l|}{scheme (i)}
$\vphantom{\int_A^A}$ \\ \hline
$\text{err}_{global}$ & $\text{ord}_{global}\!\!\!\!\!$
$\vphantom{\int_A^A}$  \\ \hline
0.166\,E$-$02 & $\vphantom{\int^A}$  \\
0.189\,E$-$05 & 6.45 \\
0.229\,E$-$07 & 6.37 \\
0.408\,E$-$09 & 5.81 \\
0.719\,E$-$11 & 5.83 \\
\hline
\end{tabular}
\caption{Error tables for the palindromic pair \method{PP 5/6 A}
         applied to problem~\eqref{eq:schroesys}. \newline
         {\em Left:}\, Local error (first step) for scheme~(i) starting with `A'
         of order~$ 5 $, and for the averaged scheme (see~\eqref{eq:pali-avaraged})
         of order~$ 6 $. \newline
         {\em Right:}\, Global error for scheme~(i) at $ t_{end} = 5.0 $.
         \label{TAB:schroesys-2}}
\end{table}

\bigskip\noindent
{\textbf{Acknowledgements.}}
{This work was supported by the Austrian Science Fund (FWF) under grant P24157-N13,
and by the Vienna Science and Technology Fund (WWTF) under grant MA-14-002.

The computational results presented have been achieved in part using the
Vienna Scientific Cluster (VSC).}


\end{document}